\documentclass[11pt, reqno]{amsart}
\usepackage{amssymb, amsthm, amsmath, amsfonts}
\usepackage[alphabetic]{amsrefs}
\usepackage{verbatim}

\newcommand{\demo}{\begin{proof}}
\newcommand{\edemo}{\end{proof}}
\newcommand{\demoname}[1]{\begin{proof}[#1]}
\newcommand{\edemoname}{\end{proof}}

\newcommand{\stepname}{Step}

\theoremstyle{plain}
\newtheorem{theorem}{Theorem}[section]

\newtheorem{conjecture}[theorem]{Conjecture}
\newtheorem{corollary}[theorem]{Corollary}
\newtheorem{lemma}[theorem]{Lemma}

\theoremstyle{definition}
\newtheorem{example}[theorem]{Example}
\newtheorem*{definition}{Definition}
\newtheorem{step}{Step}
\newtheorem{casee}{Case}
\newtheorem{stepn}{\stepname}
\newtheorem*{stepnn}{\stepname}

\newcommand{\thm}{\begin{theorem}}
\newcommand{\ethm}{\end{theorem}}
\newcommand{\conj}{\begin{conjecture}}
\newcommand{\econj}{\end{conjecture}}

\newcommand{\expl}{\begin{example}}
\newcommand{\eexpl}{\qex\end{example}}
\newcommand{\defn}{\begin{definition}}
\newcommand{\edefn}{\qef\end{definition}}
\newcommand{\defnnb}{\begin{definition}}
\newcommand{\edefnnb}{\end{definition}}
\newcommand{\stp}{\begin{step}}
\newcommand{\estp}{\end{step}}
\newcommand{\cse}{\begin{casee}}
\newcommand{\ecse}{\end{casee}}
\newcommand{\stpn}[1]{\renewcommand{\stepname}{#1}\begin{stepn}}
\newcommand{\estpn}{\end{stepn}}
\newcommand{\stpnn}[1]{\renewcommand{\stepname}{#1}\begin{stepnn}}
\newcommand{\estpnn}{\end{stepnn}}

\newcommand{\coro}{\begin{corollary}}
\newcommand{\ecoro}{\end{corollary}}
\newcommand{\lem}{\begin{lemma}}
\newcommand{\elem}{\end{lemma}}

\providecommand{\qexsymbol}{$\lozenge$}%
\newcommand{\mathqex}{\quad\hbox{\qexsymbol}}
\DeclareRobustCommand{\qex}{%
  \ifmmode \mathqex
  \else
    \leavevmode\unskip\penalty9999 \hbox{}\nobreak\hfill
    \quad\hbox{\qexsymbol}%
  \fi
}
\providecommand{\qefsymbol}{$\triangle$}%
\newcommand{\mathqef}{\quad\hbox{\qefsymbol}}
\DeclareRobustCommand{\qef}{%
  \ifmmode \mathqef
  \else
    \leavevmode\unskip\penalty9999 \hbox{}\nobreak\hfill
    \quad\hbox{\qefsymbol}%
  \fi
}

\newcommand{\enum}{\begin{enumerate}}
\newcommand{\eenum}{\end{enumerate}}

\newcommand{\nn}{\mathbb{N}}
\newcommand{\zz}{\mathbb{Z}}
\newcommand{\znn}{\nn \cup \{0\}}
\newcommand{\rrr}[1]{{\mathbb{R}}^{#1}}
\newcommand{\rr}{\mathbb{R}}

\newcommand{\nd}[1]{$#1$\nobreakdash-\hspace{0pt}}

\newcommand{\aang}{\alpha}

\newcommand{\faceof}{\prec}
\newcommand{\offace}{\succ}

\newcommand{\nsimp}[1]{\nd{#1}simplex}

\newcommand{\nface}[1]{\nd{#1}face}

\DeclareMathOperator{\starpl}{star}
\DeclareMathOperator{\linkpl}{link}
\DeclareMathOperator{\linksofl}{LI}

\newcommand{\stark}[2]{\starpl\,({#1}, {#2})}
\newcommand{\linkk}[2]{\linkpl\,({#1}, {#2})}

\newcommand{\scvf}{real-valued function}

\newcommand{\svf}{real-valued vertex-supported function}

\newcommand{\ldv}{locally determined}

\newcommand{\ldvc}{combinatorially locally determined}

\newcommand{\gld}{geometrically locally determined}

\newcommand{\simpts}{\nd{2}dimensional simplicial complexes}

\newcommand{\linksof}[1]{\linksofl({#1})}  

\newcommand{\simpcs}{simplicial complexes}

\newcommand{\phiw}[2]{\phi({#1}, {#2})}
\newcommand{\phihvkh}[3]{{#3}(\linkk {#1}{#2})}

\newcommand{\Lambdaa}[1]{\Lambda({#1})}
\newcommand{\verts}[1]{{#1}^{(0)}}

\newcommand{\mcalt}{\mathcal{T}}
\newcommand{\eeq}[1]{E_{#1}}

\begin{document}

\title{Locally Determined Functions of Finite Simplicial Complexes that are Linear Combinations of the Numbers of Simplices in Each Dimension}
\author{Ethan D.\ Bloch}
\address{Bard College\\
Annandale-on-Hudson, NY 12504\\
U.S.A.}
\email{bloch@bard.edu}
\thanks{I would like to thank Sam Hsaio for bringing this problem to my attention}
\date{}
\begin{abstract}
The Euler characteristic, thought of as a function that assigns a numerical value to every finite simplicial complex, is locally determined in both a combinatorial sense and a geometric sense.  In this note we show that not every function that assigns a numerical value to every finite simplicial complex via a linear combination of the numbers of simplices in each dimension is locally determined in either sense.  In particular, the Charney-Davis quantity $\lambda (L)$ is not locally determined in either sense if it is defined on a set of simplicial complexes that includes all flag spheres of a given odd dimension.
\end{abstract}
\maketitle

\markright{Locally Determined Functions of Finite Simplicial Complexes}

\section{Introduction}
\label{secINTR}

There are a number of contexts in which to consider the notion that a function that assigns a numerical value to every finite simplicial complex (for example the Euler characteristic) is locally determined.  We discuss two approaches here.  The common idea to both approaches is that if $\mcalt$ is a set of finite \simpcs, and if $\Lambda$ be a \scvf\ on $\mcalt$, the function $\Lambda$ is \ldv\ if there is an appropriate type of real-valued function $f$ defined at each vertex of each simplicial complex in $\mcalt$ such that
\[
\Lambdaa K = \sum_{v \in \verts K} f(v)
\]
for every $K \in \mcalt$, where $\verts K$ is the set of vertices of $K$.

In the combinatorial approach of \cite{LEVI} and \cite{FORM5}, described below in Section~\ref{secCLDF}, the number $f(v)$ depends only upon the combinatorial nature of the link of $v$ in $K$.  For example, it is observed in those two papers that the Euler characteristic on the set of all \simpcs\ is \ldvc\ by the real valued function $e$ on the set of all \simpcs\ given by 
\[
e(M) = 1 + \sum_i \frac{(-1)^{i+1}}{i+2} f_i(M)
\]
for all \simpcs\ $M$.

There is also a geometric notion of a function on a set of finite simplicial complexes being locally determined, as discussed in \cite{BL13}, and described below in Section~\ref{secGLDF}.  This  approach is inspired by the polyhedral analog of the Gauss-Bonnet Theorem, which is as follows.  Let $M$ be a finite polyhedral surface in $\rrr{3}$.  If $v \in \verts M$, the angle defect at $v$ is defined to be $d_v = 2\pi - \sum \alpha_i$, where the $\alpha_i$ are the angles of the triangles containing $v$.  The polyhedral analog of the Gauss-Bonnet Theorem says $\sum_{v \in \verts M} d_v = 2\pi \chi(M)$.  Rather than viewing this formula as stating that the angle defects at the vertices add up to something nice, it can be viewed as stating that the Euler characteristic of the surface is locally determined by a geometrically calculated quantity.  The polyhedral Gauss-Bonnet Theorem can be generalized to higher dimensions and to non-manifolds in more than one way, as seen, among others, in \cite{BA3}, \cite{C-M-S}, \cite{SH2}, \cite{GR2}, \cite{G-S2} and \cite{BL4}; we take the approach of the latter.  For Descartes' original work on the angle defect see \cite{FE}; for a very accessible treatment of the angle defect see \cite{BA1}.

In contrast to the combinatorial approach of \cite{LEVI} and \cite{FORM5}, where the local calculation at each vertex of a simplicial complex depends only upon the combinatorial nature of the link (or star) of the vertex, in the geometric approach the local calculation at a vertex depends upon local geometric information that makes use of an embedding of the simplicial complex in Euclidean space.  Hence, in the geometric approach, rather than considering simplicial complexes to be the same if they are combinatorially equivalent (as in \cite{LEVI} and \cite{FORM5}), we consider simplicial complexes that are embedded in Euclidean space, and view different embeddings of the same abstract simplicial complex as different.

The Euler characteristic is not only locally determined in both the combinatorial and geometric senses, but it is the unique function that is locally determined and that satisfies some additional conditions.  In the combinatorial approach, \cite{LEVI} and \cite{FORM5} show that the Euler characteristic is, up to a scalar multiple, the unique \ldvc\ numerical invariant of finite simplicial complexes that assigns the same number to every cone; that would hold, in particular, for a numerical invariant that is a homotopy type invariant.  In the geometric approach, \cite{BL13} shows in the \nd{2}dimensional case that the Euler characteristic is, up to a scalar multiple, the unique \gld\ numerical invariant of finite simplicial surfaces that assigns the same number to every pyramid and bipyramid.

However, whereas the Euler characteristic is a very useful example, there are combinatorially invariant ways to assign a numerical value to every finite simplicial complex that are not constant on all cones (not to mention are not homotopy type invariants).  An example of such a function is the Charney-Davis quantity $\lambda (K)$, as defined in \cite{CH-DA}.  If a function such as the Charney-Davis quantity were locally defined in either the combinatorial or the geometric approach, that might provide a useful tool for its study.  The purpose of this note, however, is to give an example that shows that not all such functions are locally determined; in particular the Charney-Davis quantity defined on odd-dimensional simplicial flag spheres is not locally determined in either sense.

We start with some notation. (For basic definitions regarding simplicial complexes, see for example \cite{HU} and \cite{R-S}.)  Throughout this note, all simplicial complexes are assumed to be finite.  Let $K$ be s finite simplicial complex.  Let $|K|$ denote the underlying space of $K$, and let $\verts K$ denote the set of vertices of $K$.  Let $f_i(K)$ denote the number of \nd{i}simplices of $K$ for each $i \in \{0, 1, \ldots, \dim K\}$; we also use the standard convention that $f_{-1}(K) = 1$, and $f_i(K) = 0$ for each $i \in \zz - \{-1, 0, 1, \ldots, \dim K\}$.  If $v \in \verts K$, let $\stark vK$ and $\linkk vK$ denote the star and link of $v$ in $K$, respectively.

The Charney-Davis quantity $\lambda (K)$ is defined to be 
\[
\lambda (K) = \sum_{i=-1}^{\dim K} \left({\textstyle -\tfrac 12}\right)^{i+1} f_i(K).
\]
Simple examples show that $\lambda$ is not constant on cones, even within a fixed dimension of simplicial complexes.  The Charney-Davis Conjecture (Conjecture~D in \cite{CH-DA}) concerns the sign of $\lambda (K)$ for certain odd-dimensional simplicial complexes, which includes all flag spheres.  We will not be referring to the Charney-Davis Conjecture here, but we do need the following definition.

\defn
Let $K$ be a simplicial complex.  The simplicial complex $K$ is a \textbf{flag complex} if for any subset $T \subseteq \verts K$, if every two distinct vertices of $T$ are joined by an edge then $T$ is the set of vertices of a face of $K$.
\edefn

The Charney-Davis quantity is an example of a real-valued function $\Lambda$ on a set of finite simplicial complexes $\mcalt$ that has the form $\Lambda (K) = \sum_{i = -1}^{\dim K} b_if_i(K)$ for all $K \in \mcalt$, for some $b_{-1}, b_0, b_1, \ldots \in \rr$.  We ask which such functions $\Lambda$ are \ldv.

Our result in the combinatorial setting is the following.

\thm\label{thmMAIN}
Let $\mcalt$ be a set of finite \simpcs, and let $\Lambda$ be the \scvf\ on $\mcalt$.  Suppose that $\Lambda$ has the form $\Lambda (K) = \sum_{i = -1}^{\dim K} b_if_i(K)$ for all $K \in \mcalt$, for some $b_{-1}, b_0, b_1, \ldots \in \rr$.
\enum
\item\label{itAC}
If all the simplicial complexes in $\mcalt$ have the same non-zero Euler characteristic, then $\Lambda$ is \ldvc.
\item\label{itAA}
If $b_{-1} = 0$, then $\Lambda$ is \ldvc.
\item\label{itAB}
If $b_{-1} \ne 0$, and if $\mcalt$ contains all flag \nd{d}spheres for some odd integer $d$ such that $d \ge 3$, then $\Lambda$ is not \ldvc. 
\eenum
\ethm

Part~(\ref{itAC}) in the above theorem would occur, for example, when the set $\mcalt$ is the set of all \nd{m}spheres for some even $m \in \nn$.

The three cases in Theorem~\ref{thmMAIN} do not exhaust all possibilities.  However, Part~(\ref{itAB}) suffices to treat the Charney-Davis quantity, as stated in Corollary~\ref{coroCD} below.

Our result in the geometric setting is the following.

\thm\label{thmMAIN2}
Let $\mcalt$ be a set of finite \simpcs, and let $\Lambda$ be the \scvf\ on $\mcalt$.  Suppose that $\Lambda$ has the form $\Lambda (K) = \sum_{i = -1}^{\dim K} b_if_i(K)$ for all $K \in \mcalt$, for some $b_{-1}, b_0, b_1, \ldots \in \rr$.
\enum
\item\label{itAA2}
If $b_{-1} = 0$, and if $\mcalt$ is a set of \nd{n}dimensional pseudomanifolds for some integer $n$ such that $n \ge 2$, then $\Lambda$ is \gld.
\item\label{itAB2}
If $b_{-1} \ne 0$, and if $\mcalt$ contains all flag \nd{d}spheres for some odd integer $d$ such that $d \ge 3$, then $\Lambda$ is not \gld. 
\eenum
\ethm

We note that in Section~2 of \cite{BL13}, it was mistakenly claimed that  the function $\Lambda$ on the set of all \simpts\ defined by $\Lambdaa K = f_2(K)$ for all $K$ is not \gld.  It is seen by Part~(\ref{itAA2}) of Theorem~\ref{thmMAIN2} that this function $\Lambda$ is \gld\ on the set of all finite \nd{2}dimensional pseudomanifolds, and, using the ideas of the remark after the proof of that theorem, it can be see that $\Lambda$ is \gld\ on the set of all \simpts.

Finally, the following corollary is an immediate result of the above two theorems. 

\coro\label{coroCD}
Let $\mcalt$ be a set of finite \simpcs\ that contains all flag \nd{d}spheres for some odd integer $d$ such that $d \ge 3$.  Then the Charney-Davis function $\lambda$ on $\mcalt$ is not \ldv\ combinatorially or geometrically. 
\ecoro

\section{Combinatorially Locally Determined Functions}
\label{secCLDF}

For the combinatorial approach of \cite{LEVI} and \cite{FORM5}, we need the following notation.  If $\mcalt$ is a set of \simpcs, let 
\[
\linksof {\mcalt} = \{\linkk vK \mid K \in \mcalt \text{ and } v \in \verts K\}.
\]

The idea of a function on a set of simplicial complexes being locally determined is that the value of the function of a simplicial complex equals the sum of the values of some other function calculated in a ``neighborhood'' of each vertex of the simplicial complex.  The standard notion of neighborhood of a vertex in a simplicial complex is the star of the vertex, but the star of the vertex is determined by the link of the vertex, and in the following definition, which is from \cite{LEVI} and \cite{FORM5}, it is convenient to use the link rather than the star. 

\defn
Let $\mcalt$ be a set of \simpcs, and let $\Lambda$ be a \scvf\ on $\mcalt$.  The function $\Lambda$ is \textbf{\ldvc} by a \scvf\ $h$ on $\linksof {\mcalt}$ if $h$ is invariant under combinatorial equivalence and if
\[
\Lambdaa K = \sum_{v \in \verts K} h(\linkk vK)
\]
for every $K \in \mcalt$.
\edefn

The adjective ``combinatorially'' in the above definition is not used in \cite{LEVI} and \cite{FORM5}, but for the sake of clarity it seems appropriate to use it at present.

For the following proof, we need the following basic facts about joins of simplicial complexes.  Let $K$ and $L$ be finite simplicial complexes.  Let $K \ast L$ denote the join of $K$ and $L$.  Then $\dim (K \ast L) = \dim K + \dim L + 1$, and  
\begin{equation}\label{eqABA}
f_r(K \ast L) = \sum_{i=-1}^{r} f_{r-i-1}(K)f_i(L).
\end{equation}
for each $r \in \zz$.  Additionally, if $v \in \verts K$, then 
\begin{equation}\label{eqABO}
\linkk {v}{K \ast L} = \linkk vK \ast L.
\end{equation}

\demoname{Proof of Theorem~\ref{thmMAIN}}
Parts~(\ref{itAC}) and (\ref{itAA}) are very simple.  For each case, we will define a sequence $a_{-1}, a_0, a_1, \ldots \in \rr$, and we will define a \scvf\ $g$ on $\linksof {\mcalt}$ that has the form $g(M) = \sum_{i = -1}^{\dim M} a_if_i(M)$ for each $M \in \linksof {\mcalt}$.

First, we make the following observation.  Let $K \in \mcalt$, let $v \in \verts K$ and let $i \in \{-1, \ldots, \dim K - 1\}$.  It is  straightforward to see that $\sum_{v \in \verts K} f_i(\linkk vK) = (i + 2)f_{i+1}(K)$.  Then
\begin{equation}\label{eqEBG}
\begin{aligned}
\sum_{v \in \verts K} \phihvkh vKg &= \sum_{v \in \verts K} \sum_{i = -1}^{\dim K - 1} a_if_i(\linkk vK)\\
 &= \sum_{i = -1}^{\dim K - 1} a_i\sum_{v \in \verts K} f_i(\linkk vK)\\
 & = \sum_{i = -1}^{\dim K - 1} a_i(i + 2)f_{i+1}(K).
\end{aligned}
\end{equation}

For Part~(\ref{itAC}), where we assume that all the simplicial complexes in $\mcalt$ have Euler characteristic $E$, for some $E \ne 0$, let $a_i = \frac {b_{i+1}}{i+2} + (-1)^{i+1}\frac {b_{-1}}{E(i+2)}$ for all $i \in \{-1, 0, 1, \ldots\}$.  For Part~(\ref{itAA}), where we assume that  that $b_{-1} = 0$, let $a_i = \frac {b_{i+1}}{i+2}$ for all $i \in \{-1, 0, 1, \ldots\}$.  In both cases, it is straightforward to verify that Equation~(\ref{eqEBG}) implies  $\sum_{v \in \verts K} \phihvkh vKg = \sum_{i = -1}^{\dim K} b_if_i(K) = \Lambda (K)$; the details are omitted.

For Part~(\ref{itAB}), suppose that $b_{-1} \ne 0$, and that there is some odd integer $d$ such that $d \ge 3$ and that $\mcalt$ contains all flag \nd{d}spheres.  Suppose further that $\Lambda$ is \ldv\ by a \scvf\ $h$ on $\linksof {\mcalt}$.

Let $n, m \in \nn$ be such that $n, m \ge 4$ and $n \ne m$.  Let $C_n$ denote the cycle with $n$ vertices (that is, the graph).

Let $s, t \in \znn$.  Let
\[
T_{s,t} = \underbrace{C_n \ast \cdots \ast C_n}_{s \text{ times}} \ast \underbrace{C_m \ast \cdots \ast C_m}_{t \text{ times}}.
\]
Then $T_{s,t}$ is a \nd{(2(s+t) - 1)}dimensional simplicial complex.  We note that $T_{s,t}$ is a flag complex because $C_n$ and $C_m$ are flag complexes (because $n, m \ge 4$), and the join of flag complexes is a flag complex (Item~2.7.1 of \cite{CH-DA}).  We also see that $T_{s,t}$ is a sphere, because the join of spheres is a sphere (Lemma~1.13 of \cite{HU}).

It can be verified that
\[
f_i(T_{s,0}) = \sum_{j=0}^s \tbinom sj \tbinom j{2j-i-1} n^j
\]
for all $i \in \zz$; if $i \in \{-1, 0, \ldots, 2s - 1\}$, count the ways \nd{i}dimensional simplices in this join can be formed, and if $i \notin \{-1, 0, \ldots, 2s - 1\}$, then the formula correctly yields $f_i(T_{s,0}) = 0$.  A similar formula holds for $f_i(T_{0,t})$, and then, using Equation~(\ref{eqABA}), we see that 
\begin{equation}\label{eqABD}
\begin{aligned}
f_r(T_{s,t}) &= f_r(T_{s,0} \ast T_{0,t}) = \sum_{i=-1}^{r} f_{r-i-1}(T_{s,0})f_i(T_{0,t})\\
 &= \sum_{i=-1}^{r}\sum_{j=0}^s\sum_{k=0}^t \tbinom sj \tbinom j{2j+i-r} \tbinom tk \tbinom k{2k-i-1} n^jm^k
\end{aligned}
\end{equation}
for all $r \in \zz$.

Let $v \in \verts {(T_{s,t})}$.  Then $v$ is either in a copy of $C_n$ or a copy of $C_m$.  Let $S^0$ denote a two-element set.  If $v$ is in a copy of $C_n$, then by Equation~(\ref{eqABO}) we see that $\linkk v{T_{s,t}}$ is isomorphic to $S^0 \ast T_{s-1,t}$, and we let $A_{s-1,t} = \phihvkh v{T_{s,t}}h = h(S^0 \ast T_{s-1,t})$; similarly, if $v$ is in a copy of $C_m$, we let $A_{s,t-1} = \phihvkh v{T_{s,t}}h = h(S^0 \ast T_{s,t-1})$.

Let $p \in \nn$ be such that $d = 2p - 1$.  Then $p \ge 2$.  Let $u \in \{0, 1, \ldots, p\}$.  Then $T_{p-u,u}$ is a \nd{(2p-1)}dimensional simplicial flag sphere, and hence $T_{p-u,u} \in \mcalt$.  Using the fact that $\Lambda$ is \ldv\ by $h$, we have
\begin{equation}\label{eqABB}
\begin{aligned}
\sum_{r = -1}^{2p-1} b_rf_r(T_{p-u,u}) &= \Lambda (T_{p-u,u})  = \sum_{v \in \verts {(T_{p-u,u})}} \phihvkh v{T_{p-u,u}}h\\
 &= f_0(T_{p-u,0})A_{p-u-1,u} + f_0(T_{0,u})A_{p-u,u-1}.
\end{aligned}\tag{$\eeq u$}
\end{equation}

Next, we take Equations~($\eeq 0$), ($\eeq 1$), \ldots, ($\eeq p$), and form the linear combination $\sum_{w=0}^{p} (-1)^w\tbinom pw (\frac mn)^{p-w} \eeq w$, which, after rearranging, yields
\begin{equation}\label{eqABC}
\begin{aligned}
\sum_{r = -1}^{2p-1} b_r &\sum_{w=0}^{p} (-1)^w \tbinom pw \left(\tfrac mn\right)^{p-w}f_r(T_{p-w,w})\\
 &= \sum_{w=0}^{p} (-1)^w\tbinom pw \left(\tfrac mn\right)^{p-w}\left[f_0(T_{p-w,0})A_{p-w-1,w} + f_0(T_{0,w})A_{p-w,w-1}\right].
\end{aligned}
\end{equation}

We now simplify Equation~(\ref{eqABC}), starting with the right-hand side of the equation, which is a linear combination of terms of the form $A_{p-u,u-1}$, where $u \in \{0, 1, \ldots, p + 1\}$.  Each such term appears twice, once in each of Equations~($\eeq {u-1}$) and ($\eeq u$).  (It might be thought that each of $A_{p,-1}$ and $A_{-1,p}$ appear only once, in Equations~($\eeq {0}$) and ($\eeq p$), respectively, but each of $A_{p,-1}$ and $A_{-1,p}$ has coefficient $f_0(T_{0,0})$, which is zero, and so we can ignore these terms.)  The sum of the two coefficients of $A_{p-u,u-1}$ from Equations~($\eeq {u-1}$) and ($\eeq u$) is
\begin{align*}
(-1)^{u-1}&\tbinom p{u-1} \left(\tfrac mn\right)^{p-(u-1)}f_0(T_{p-(u-1),0}) + (-1)^u\tbinom pu \left(\tfrac mn\right)^{p-u}f_0(T_{0,u})\\
 &=(-1)^{u-1}\left(\tfrac mn\right)^{p-u}\left[\tbinom p{u-1} \left(\tfrac mn\right)(p-u+1)n - \tbinom pu um\right]\\
 &=0,
\end{align*}
where the last equality can be verified easily.  We deduce that the right-hand side of Equation~(\ref{eqABC}) is zero.

Next, let $r \in \{-1, 0, \ldots, 2p-1\}$, and let $G_r$ denote the coefficient of $b_r$ in the left-hand side of Equation~(\ref{eqABC}).

Because $f_{-1}(T_{p-w,w}) = 1$ for all $w \in \{0, 1, \ldots, p\}$, we see that
\[
G_{-1} = \sum_{w=0}^{p} (-1)^w \tbinom pw \left(\tfrac mn\right)^{p-w} = \left(\tfrac mn - 1\right)^p.
\]

Now suppose that $r \ne -1$.  Then, using Equation~(\ref{eqABD}), we see that
\begin{align}
G_r &= \sum_{w=0}^{p} (-1)^w \tbinom pw \left(\tfrac mn\right)^{p-w}f_r(T_{p-w,w})\notag\\
 &= \sum_{w=0}^{p} (-1)^w \tbinom pw \left(\tfrac mn\right)^{p-w} \sum_{i=-1}^{r}\sum_{j=0}^{p-w}\sum_{k=0}^w \tbinom {p-w}j \tbinom j{2j+i-r} \tbinom wk \tbinom k{2k-i-1} n^jm^k.\notag\\
 &= \sum_{w=0}^{p} \sum_{i=-1}^{r}\sum_{j=0}^{p-w}\sum_{k=0}^w (-1)^w \tbinom pw \tbinom {p-w}j \tbinom j{2j+i-r} \tbinom wk \tbinom k{2k-i-1} \tfrac {m^{p - w + k}}{n^{p - w - j}}\label{eqABE}.
\end{align}

Observe that $G_r$ is a Laurent polynomial in $m$ and $n$.  We will use the substitution $a = p - w + k$ and $b = p - w - j$.  Because $0 \le k \le w$ and $0 \le j \le p - w$, then $0 \le a \le p$ and $0 \le b \le p$; it also follows that $a - (p - w) = k \ge 0$ and $(p - w) - b = j \ge 0$.  It is evident that $a \ge b$.  In fact, there is no term in $G_r$ that has $\frac {m^a}{n^b}$ with $a = b$.  Suppose the contrary.  Because $k, j \ge 0$, the only values of $k$ and $j$ that would yield $a = b$ are $k = 0 = j$.  If that were the case, then $\tbinom j{2j+i-r} \tbinom k{2k-i-1} = \tbinom 0{i-r} \tbinom 0{-i-1}$.  Note that $i \in \{-1, 0, \ldots, r\}$.  If $i < r$, then $\tbinom 0{i-r} = 0$, and if $i > -1$, then $\tbinom 0{-i-1} = 0$.  Hence, the coefficients of $\frac {m^a}{n^b}$ with $a = b$ are all zero, and we may therefore restrict our attention to the case where $a > b$.  

Let $D_{a,b}$ denote the coefficient of $\frac {m^a}{n^b}$ in Equation~(\ref{eqABE}).  For each possible value of $w$, there is one choice of each of $k$ and $j$ that yield the desired powers of $m$ and $n$.  Specifically, it is seen that
\begin{equation}\label{eqABG}
D_{a,b} = \sum_{w=0}^{p} (-1)^w \tbinom pw \tbinom {p-w}{(p-w)-b} \tbinom w{a-(p-w)}\sum_{i=-1}^{r} \tbinom {(p-w)-b}{2(p-w)-2b+i-r} \tbinom {a-(p-w)}{2a-2(p-w)-i-1}.
\end{equation}

Let $w \in \{0, 1, \ldots, p\}$.  Using the definition of binomial coefficients, it is easy to verify that
\begin{equation}\label{eqABL}
\tbinom pw \tbinom {p-w}{(p-w)-b} \tbinom w{a-(p-w)} = \tbinom p{a-b} \tbinom {p-(a-b)}{b} \tbinom {a-b}{a-(p-w)}.
\end{equation}

Next, recalling that $a - (p - w) \ge 0$ and $(p - w) - b \ge 0$, we observe that if $i \in \zz$, then $i < -1$ implies $\tbinom {a-(p-w)}{2a-2(p-w)-i-1} = 0$, and $i > r$ implies $\tbinom {(p-w)-b}{2(p-w)-2b+i-r} = 0$.  Let $c, e \in \zz$, and suppose that $c \le -1$ and $e \ge r$.  It then follows that that 
\begin{equation}\label{eqABH}
\sum_{i=-1}^{r} \tbinom {(p-w)-b}{2(p-w)-2b+i-r} \tbinom {a-(p-w)}{2a-2(p-w)-i-1} = \sum_{i=c}^{e} \tbinom {(p-w)-b}{2(p-w)-2b+i-r} \tbinom {a-(p-w)}{2a-2(p-w)-i-1}.
\end{equation}

Let $A = 2(p-w) - 2b - r$ and $B = 2a - 2(p-w) - 1$.  Then 
\[
\sum_{i=c}^{e} \tbinom {(p-w)-b}{2(p-w)-2b+i-r} \tbinom {a-(p-w)}{2a-2(p-w)-i-1} = \sum_{v=A+c}^{A+e} \tbinom {(p-w)-b}{v} \tbinom {a-(p-w)}{(A + B) - v}.
\]
If we choose values of $c$ and $e$ so that $\{0, 1, \ldots, A+B\} \subseteq \{A+c, A+c+1, \ldots, A+e\}$, we can apply Vandermonde's Convolution Formula to deduce that 
\begin{equation}\label{eqABJ}
\sum_{i=-1}^{r} \tbinom {(p-w)-b}{2(p-w)-2b+i-r} \tbinom {a-(p-w)}{2a-2(p-w)-i-1} = \tbinom {a-b}{A+B} = \tbinom {a-b}{2(a-b)-r-1}.
\end{equation}

Combining Equations~(\ref{eqABG}), (\ref{eqABL}) and (\ref{eqABJ}), and using the substitution $z = a  - (p - w)$, yields
\begin{align*}
D_{a,b} &= \sum_{w=0}^{p} (-1)^w \tbinom p{a-b} \tbinom {p-(a-b)}{b} \tbinom {a-b}{a-(p-w)}\tbinom {a-b}{2(a-b)-r-1}\\
 &= \tbinom p{a-b} \tbinom {p-(a-b)}{b} \tbinom {a-b}{2(a-b)-r-1} \sum_{z=a-p}^{a} (-1)^{p - (a - z)} \tbinom {a-b}{z},
\end{align*}
Observe that $a - p \le 0 < a - b \le a$.  Hence, the only values of $z$ for which $\tbinom {a-b}{z}$ is non-zero are $z \in \{0, \ldots, a - b\}$.  Then
\begin{align*}
D_{a,b} &= \tbinom p{a-b} \tbinom {p-(a-b)}{b} \tbinom {a-b}{2(a-b)-r-1} (-1)^{p - a} \sum_{z=0}^{a-b} (-1)^{z} \tbinom {a-b}{z}\\
 &= \tbinom p{a-b} \tbinom {p-(a-b)}{b} \tbinom {a-b}{2(a-b)-r-1} (-1)^{p - a} (1-1)^{a-b} = 0.
\end{align*}

Putting the above together, we see that Equation~(\ref{eqABC}) reduces to the very simple equation $\left(\frac mn - 1\right)^pb_{-1} = 0$.  Given that $m \ne n$ and $b_{-1} \ne 0$, we have reached a contradiction, from which we deduce that $\Lambda$ is not \ldv.
\edemoname

We note that whereas in Part~(\ref{itAB}) of Theorem~\ref{thmMAIN} it is hypothesized that $\mcalt$ contains all flag \nd{d}spheres for some odd integer $d$ such that $d \ge 3$, it is seen in the proof of the theorem that not all flag \nd{d}spheres are needed, but rather, by using $n = 4$ and $m = 5$, it would suffice to include only those flag \nd{d}spheres with up to $\frac {5(d+1)}2$ vertices.  The proof of the theorem was given with arbitrary $n$ and $m$, rather than only $n = 4$ and $m = 5$, because it is easier to see what is going on by treating the more general case.

\section{Geometrically Locally Determined Functions}
\label{secGLDF}

For the geometric approach of \cite{BL13}, we need the following definitions (which, in contrast to the original, are given here for arbitrary dimensions).  Recall that in this approach, we consider simplicial complexes that are embedded in Euclidean space, and view different embeddings of the same abstract simplicial complex as different simplicial complexes.

\defn 
Let $\mcalt$ be a set of simplicial complexes.  A \textbf{\svf} on $\mcalt$ is a function $\phi$ that assigns to every $K \in \mcalt$, and every $v \in \verts K$, a real number $\phiw vK$.
\edefn

For the following definition, suppose that $K$ and $\displaystyle \{K_n\}_{n=1}^\infty$ are combinatorially equivalent simplicial complexes, and all are embedded in the same Euclidean space.  We can think of these simplicial complexes as embeddings of the same abstract simplicial complex.  We write $\lim_{n\to\infty} K_n = K$ to denote pointwise convergence of these embeddings; it suffices to verify convergence at the vertices of the abstract simplicial complex.

\defnnb
Let $\mcalt$ be a set of simplicial complexes, and let $\phi$ be a \svf\ on $\mcalt$.  
\enum
\item
The function $\phi$ is \textbf{invariant under subdivision} if the following condition holds.  If $K, J \in \mcalt$, where $J$ is a subdivision of $K$, and if $v \in \verts K$, then $\phiw vK = \phiw vJ$.
\item
The function $\phi$ is \textbf{invariant under simplicial isometries of stars} if the following condition holds.  If $K, L \in \mcalt$, if $v \in \verts K$ and $w \in \verts L$, and if there is a simplcial isometry $|\stark vK| \to |\stark wL|$ that takes $v$ to $w$, then $\phiw vK = \phiw wL$.
\item
The function $\phi$ is \textbf{continuous} if the following condition holds.  Let $K$ and $\displaystyle \{K_n\}_{n=1}^\infty$ be combinatorially equivalent simplicial complexes in $\mcalt$, all embedded in the same Euclidean space.  Suppose $\lim_{n\to\infty} K_n = K$.  If $v \in \verts K$, and if the corresponding verticex of $K_n$ is labeled $v_n$, then $\lim_{n\to\infty}\phiw {v_n}{K_n} = \phiw vK$.\qef  
\eenum
\edefnnb

\defn
Let $\mcalt$ be a set of \simpcs, and let $\Lambda$ be a \scvf\ on $\mcalt$.  The function $\Lambda$ is \textbf{\gld} by a \svf\ $\phi$ on $\mcalt$ if $\phi$ is invariant under simplicial isometries of stars, is invariant under subdivision and is continuous, and if
\[
\Lambdaa K = \sum_{v \in \verts K} \phiw vK
\]
for every $K \in \mcalt$.
\edefn

For Part~(\ref{itAA2}) of the following proof, which is a simple variation of an argument in \cite{BL4}, we adopt the convention that all angles are normalized so that the volume of the unit \nd{(n-1)}sphere in \nd{(n-1)}measure is $1$ in all dimensions.  For any \nsimp n\ $\sigma^n$ in Euclidean space, and any \nface i\ $\eta^i$ of $\sigma^n$, let $\aang(\eta^i, \sigma^n)$ denote the solid angle in $\sigma^n$ along $\eta^i$, where by normalization such an angle is a number in $[0, 1]$.  

We make use here of a lemma, found in many places and stated as Lemma~3.1 in \cite{BL4}, which generalizes the fact that the angles of a planar triangle add up to $\pi$ (or $1/2$ when normalized); the lemma reduces to that result when $n = 2$.  A simple proof of this lemma appears on p.~24 of \cite{HF}; for historical remarks, see p.~312 of \cite{GR1}.

\demoname{Proof of Theorem~\ref{thmMAIN2}}
For Part~(\ref{itAA2}), let $K \in \mcalt$.  Because $K$ is an \nd{n}dimensional pseudomanifold, we have $(n + 1)f_n(K) = 2 f_{n-1}(K)$.

We use the notation $\eta^i$ to denote an \nsimp i\ of $K$.  Let 
\[
P = \frac{2(-1)^n}{n - 1} \left[\frac {(n + 1)}{2}b_{n-1} + b_n\right].
\]

For each $v \in \verts K$, let 
\[
\phiw vK = \sum_{i=0}^{n-2} 
\frac{1}{i+1}   \sum_{\eta^i \ni v} \left[b_i + (-1)^i P \sum_{\sigma^n \offace \eta^i} 
\aang(\eta^i, \sigma^n)\right].
\]
Because $\phi$ is determined by angle sums of the form $\sum_{\sigma^n \offace \eta^i} \aang(\eta^i, \sigma^n)$, then it is invariant under simplicial isometries of stars, it is invariant under subdivision and it is continuous.

We compute
\allowdisplaybreaks
\begin{align*}
\sum_{v \in \verts K} \phiw vK &= \sum_{v \in \verts K} \sum_{i=0}^{n-2} 
\frac{1}{i+1}   \sum_{\eta^i \ni v} \left[b_i + (-1)^i P \sum_{\sigma^n \offace \eta^i} 
\aang(\eta^i, \sigma^n)\right]\\
 &= \sum_{i=0}^{n-2} \frac{1}{i+1} \sum_{\eta^i \in K} 
\sum_{v \in \eta^i} \left[b_i + (-1)^i P \sum_{\sigma^n \offace \eta^i} 
\aang(\eta^i, \sigma^n)\right]\\
 &= \sum_{i=0}^{n-2} \sum_{\eta^i \in K}  
\left[b_i + (-1)^i P \sum_{\sigma^n \offace \eta^i} 
\aang(\eta^i, \sigma^n)\right]\\
\intertext{\hfill because $\eta^i$ has $i+1$ vertices }
 &= \sum_{i=0}^{n-2} b_if_i(K) + 
P \sum_{i=0}^{n-2} \sum_{\eta^i \in K} \sum_{\sigma^n \offace \eta^i} (-1)^i \aang(\eta^i, \sigma^n) \\
 &= \sum_{i=0}^{n-2} b_if_i(K) + 
P \sum_{\sigma^n \in K} \sum_{\substack{\eta^i \faceof \sigma^n \\ 0 \leq i \leq n-2}} (-1)^i \aang(\eta^i, \sigma^n) \\
 &= \sum_{i=0}^{n-2} b_if_i(K) + 
P \sum_{\sigma^n \in K} \frac{(-1)^n(n-1)}2\\ 
\intertext{\hfill by Lemma~3.1 of \cite{BL4}}
 &= \sum_{i=0}^{n-2} b_if_i(K) + P \frac{(-1)^n (n-1)}{2} f_n(K)\\
 &= \sum_{i=0}^{n-2} b_if_i(K) + \frac {(n + 1)}{2}b_{n-1} f_n(K) + b_n f_n(K)\\
 &= \sum_{i = 0}^{n} b_if_i(K).
\end{align*}
Hence $\Lambda$ is \gld\ by $\phi$.

For Part~(\ref{itAB2}), we simply need to modify the proof of Part~(\ref{itAB}) of Theorem~\ref{thmMAIN} very slightly, as follows.  First, embed each copy of $C_n$ or $C_m$ in $\rrr{2}$ by having the vertices be on the unit circle in $\rrr{2}$, equally spaced, and then construct $T_{p-u,u}$ in $(\rrr{2})^p = \rrr{2p}$.  It is then seen that all the vertices in $T_{p-u,u}$ that are in copies of $C_n$ have isometric stars, and similarly for all the vertices in $T_{p-u,u}$ that are in copies of $C_n$.  Hence, if $\Lambda$ were \gld\ by a \svf\ $\phi$, and if $\phi$ is assumed to be invariant under simplicial isometries of stars (it does not even have to satisfy the other two conditions in the definition of \gld), then the same example used in the proof of Part~(\ref{itAB}) of Theorem~\ref{thmMAIN}, but with $\phihvkh v{T_{s,t}}h$ replaced by $\phiw v{T_{s,t}}$, will yield the exact same contradiction.
\edemoname

Finally, we note that whereas our approach in the proof of Part~(\ref{itAA2}) of Theorem~\ref{thmMAIN2} worked for \nd{n}dimensional pseudomanifolds but not all simplicial complexes, it is possible to modify the definition of $\phiw vK$ in the proof in such a way that it works for all pure finite \nd{n}dimensional simplicial complexes, though at the cost that instead of obtaining expressions of the form $\sum_{i = 0}^{n} b_if_i(K)$, each number $f_i(K)$ would be replaced by a variant of it that is weighted by the extent to which $K$ is not a pseudomanifold, using the methodology of \cite{BL4}.  We omit the details.


\begin{bibdiv}

\begin{biblist}[\normalsize]

\bib{BA1}{article}{
author = {Banchoff, Thomas},
title = {Critical points and curvature for embedded polyhedra},
journal = {J. Diff. Geom.},
volume = {1},
date = {1967},
pages = {245--256}
}

\bib{BA3}{article}{
author = {Banchoff, Thomas},
title = {Critical points and curvature for embedded polyhedra, II},
journal = {Progress in Math.},
volume = {32},
date = {1983},
pages = {34--55}
}

\bib{BL4}{article}{
author = {Bloch, Ethan D.},
title = {The angle defect for arbitrary polyhedra},
journal = {Beitr{\"a}ge Algebra Geom.},
volume = {39},
date = {1998},
pages = {379--393}
}

\bib{BL13}{article}{
author = {Bloch, Ethan D.},
title = {A characterization of the angle defect and the {E}uler characteristic in dimension $2$},
journal = {Discrete Comput. Geom.},
volume = {43},
date = {2010},
pages = {100--120}
}

\bib{CH-DA}{article}{
   author={Charney, Ruth},
   author={Davis, Michael},
   title={The Euler characteristic of a nonpositively curved, piecewise
   Euclidean manifold},
   journal={Pacific J. Math.},
   volume={171},
   date={1995},
   number={1},
   pages={117--137}
}

\bib{C-M-S}{article}{
author = {Cheeger, J.},
author = {Muller, W.},
author = {Schrader, R.},
title = {On the curvature of piecewise flat spaces},
journal = {Commun. Math. Phys.},
volume = {92},
date = {1984},
pages = {405--454}
}

\bib{FE}{book}{
author = {Federico, P. J.},
title = {Descartes on Polyhedra},
publisher = {Springer-Verlag},
address = {New York},
date = {1982}
}

\bib{FORM5}{article}{
author = {Forman, Robin},
title = {The {E}uler characteristic is the unique locally determined numerical invariant of finite simplicial complexes which assigns the same number to every cone},
journal = {Discrete Comput. Geom.},
volume = {23},
number = {4},
date = {2000},
pages = {485--488}
}

\bib{GR1}{book}{
author = {Gr{\"u}nbaum, Branko},
title = {Convex Polytopes},
publisher = {John Wiley \& Sons},
address = {New York},
date = {1967}
}

\bib{GR2}{article}{
author = {Gr{\"u}nbaum, Branko},
title = {{G}rassman angles of convex polytopes},
journal = {Acta. Math.},
volume = {121},
date = {1968},
pages = {293--302}
}

\bib{G-S2}{article}{
author = {Gr\"unbaum, Branko},
author = {Shephard, G. C.},
title = {{D}escartes' theorem in $n$ dimensions},
journal = {Enseign. Math. (2)},
volume = {37},
date = {1991},
pages = {11--15}
}

\bib{HF}{book}{
author = {Hopf, Heinz},
title = {Differential Geometry in the Large},
series = {Lecture Notes in Math.},
number = {1000},
publisher = {Springer-Verlag},
address = {Berlin},
date = {1983}
}

\bib{HU}{book}{
author = {Hudson, J. F. P.},
title = {Piecewise Linear Topology},
publisher = {Benjamin},
address = {Menlo Park, CA},
date = {1969}
}

\bib{LEVI}{article}{
author = {Levitt, Norman},
title = {The {E}uler characteristic is the unique locally determined 
numerical homotopy invariant of finite complexes},
journal = {Discrete Comput. Geom.},
volume = {7},
date = {1992},
pages = {59--67}
}

\bib{R-S}{book}{
author = {Rourke, C.},
author = {Sanderson, B.},
title = {Introduction to Piecewise-Linear Topology},
series = {Ergebnesse der Mathematik},
number = {69},
publisher = {Springer-Verlag},
address = {New York},
date = {1972}
}

\bib{SH2}{article}{
author = {Shephard, G. C.},
title = {Angle deficiencies of convex polytopes},
journal = {J. London Math. Soc.},
volume = {43},
date = {1968},
pages = {325--336}
}


\end{biblist}

\end{bibdiv}
\end{document}